# CALCULATING COST DISTRIBUTIONS OF A MULTISERVICE LOSS SYSTEM

Dr. Jorma Jormakka[1], Mr. Sourangshu Ghosh[2]

**Abstract:** Congestion pricing has received lots of attention in the scientific discussion. Congestion pricing means that the operator increases prices at the time of congestion and the traffic demand is expected to decrease. In a certain sense, shadow prices are an optimal way of congestion pricing: users are charged shadow prices, i.e., the expectations of future losses because of blocked connections. The shadow prices can be calculated exactly from Howard's equation, but this method is difficult. The paper presents simple approximations to the solution of Howard's equation and a way to derive more exact approximations. If users do not react by lowering their demand, they will receive higher bills to pay. Many users do not react to increased prices but would want to know how the congestion pricing mechanism affects the bills. The distribution of the price of a connection follows from knowing the shadow prices and the probability of a congestion state. There is another interesting distribution. The network produces profit to the operator, or equivalently, blocked connections produce a cost to the operator. The average cost rate can be calculated from Howard's equation, but the costs have some distribution. The distribution gives the risk that the actual costs exceed the average costs, and the operator should include this risk to the prices. The main result of this paper shows how to calculate the distribution of the costs in the future for congestion pricing by shadow prices and for congestion pricing with a more simple pricing scheme that produces the same average costs.

## I. INTRODUCTION

CONGESTION pricing is a state-dependent pricing strategy where the prices are higher when the network is congested. Congestion pricing has been studied extensively as a congestion control mechanism for elastic traffic demand: if users react to increased prices by lowering their demand, then the operator can use prices as a congestion control mechanism. Most authors consider congestion pricing of elastic traffic using the Game Theory [1], [2], [3], [4], and [5]. Shadow prices turn out to be the congestion prices that optimize the utility functions of the operator and the users.

The main focus in this paper is in distributions of costs when congestion pricing is applied. For a business user the interesting question is the cost of a connection. The bills of the business users are followed up and the company tries to notice errors in billing by the operator, abuse of the connections in some way, and malfunctioning of equipment, all of which can be the reason for exceptionally high bills. If a user comes to the network at a random time without knowing the congestion stage of the network and if the operator uses congestion pricing, then he will receive variable bills for the same service. Thus, congestion pricing is one cause to variability of the bill and both the company and the operator should know the effect of this mechanism. The expected distribution of the cost of a connection is obtained from the average relative costs and the state probability. Average relative costs have been calculated by many authors, the closest to this work are the calculations in [7], [10], [11] and [12]. The presented paper gives a simple approximate solution that does not require matching of parameters and shows the shape of the cost function better. A series form solution for deriving better approximation is also given.

It is important for an operator using shadow prices to know how much variability there can be in the costs over a longer time period. This is not the same as variability of the shadow prices. It is a question of the limit cost distribution. The main contribution of the paper is calculation of the cost distribution. Cost distributions have not been calculated before in the literature to Author's knowledge. Calculating the cost distribution for shadow prices turns out to be difficult, but a simpler pricing scheme can be created. For a user and the operator this simpler scheme gives practically the same prices and costs as shadow prices, but the cost distribution can be expressed by a closed formula.

We will not assume that the majority of users respond to congestion prices by lowering their traffic demand. There are reasons for dropping this assumption. Firstly, the users must notice that prices have been increased. Many users are a bit lazy in reacting to information of this type, so probably there should be an automatic agent in the user's terminal but that is unlikely to happen in the near future. Secondly, the user should want to reduce his demand. If there is an alternative network available the user may try to switch to the alternative connection. Often the user is connected to one operator's network only and cannot easily respond to price increases by shifting traffic to another network. In such a situation he should reduce his traffic demand. The need for congestion control is mainly from real-time services, today voice calls and in the future multimedia calls. The peak traffic of these services has traditionally been in the business days in the busy periods. There is no reason to expect that this situation would essentially change in the future. Bills of business users are paid by the company and both the users and the company are in the opinion that business time is so valuable that operator bills should not interfere with the ability to communicate when needed. It is unlikely that business users will respond to increased prices by lowering their demand. Some other users may indeed lower their demand and call blocking may be reduced or quality of service is improved by this method – therefore an operator may want to use congestion pricing. For a business user the main visible effect of congestion pricing is the variability of the cost of a connection.

Following [6] we treat the system as a multiservice loss system since the future applications where congestion pricing could be appropriate are most probably multimedia calls or similar services. While packet networks are queuing systems, many services still require sufficient effective bandwidth and treating the network as a loss system often gives a fairly good approximation also for packet networks. Thus, we will not study congestion pricing as a congestion control mechanism of elastic traffic in a packet network. The setting is more similar to that in [6], where possibly state-dependent prices are used in a multiservice loss network, or [7] presenting a method to calculate shadow prices a multiservice loss network.

The paper is organized as follows. In Section II multiservice loss system is described. Section III A. shows how to calculate average costs using shadow prices. This calculation also gives the distribution of costs of connections in a user's bill. Section III B. presents cost distributions to an operator over a longer time if shadow prices are strictly followed. The results are complicated and solution can be given only to one simple blocking policy. Section IV has the main result of the paper: an exact solution to the cost distribution for a charging policy, which gives the same average costs as shadow prices and for a user is practically identical.

## II. MULTISERVICE LOSS SYSTEM

A multiservice loss system (see e.g [8]) has a number of different connection types indexed by $j \in \{1,...,K\}$. Each has own bandwidth demand: a call of type $j$ needs $b_j$ units of capacity. Each call type has a Poisson arrival process with an arrival rate $\lambda_j$. The call holding time for each call type is negatively exponentially distributed with the service rate $\mu_j$. Let the set $\Omega$ be the set of admitted state vectors $\bar{q}$. The policy sets $R_{\bar{q}} = \{j | \bar{q} + e_j \in \Omega\}$ define whether a call is accepted. Here $\bar{q} + e_j = (q_1,...q_j + 1,...q_K)$ where $e_j$ has 1 in the index $j$, i.e. $e_j = (0,\ldots,1,\ldots)$. The policy "a new call of type $j$ is accepted in the state $\bar{q}$" is expressed using policy sets with the condition $j \in R_{\bar{q}}$.

Example 1. A call of type $k$ takes $d_k$ capacity units and the maximum number of capacity units is $C$. If the total number of used capacity units

$$c = \sum_{k=1}^{K} d_k q_k = \bar{q}^T \bar{d} \text{ is larger than } C - d_k, \quad (1)$$

a call of type $k$ is blocked. Otherwise it is accepted. Thus

$$R_{\bar{q}} = \{j | \bar{q}^T \bar{d} \le C - d_j\}.$$

Example 2: A call of type $k$ takes $d_k$ capacity units. If the number of used capacity units by calls of type $k$ exceeds a threshold $C_k$ then a call of type $k$ is blocked.

The policy in Example 2 is less efficient than the policy in Example 1 because not all capacity is available to all call types.

The instantaneous transition matrix for the multi-service loss system can be written down by looking at the system at the state $\bar{q}=(q_1,...,q_K)$ on a moment of time *n*. Let the total time be *T* and let there be *N* discrete time steps in time *T*. (We let $N\to\infty$.) Then we have the following transitions:

- a connection of type $j$ is accepted $(n,\bar{q})\to(n+1,\bar{q}+e_j)$ with probability $\frac{T}{N}\lambda_j 1_{j\in R_{\bar{q}}}$ ,

- the state does not change $(n,\bar{q})\to(n+1,\bar{q})$ with probability

$$1-\frac{T}{N}\sum_{j\in R_{\bar{q}}}\lambda_j-\frac{T}{N}\sum_{j=1}^{K}\mu_j q_j ,$$

- a connection of type $j$ finishes $(n,\bar{q})\to(n+1,\bar{q}-e_j)$ with probability $\frac{T}{N}\mu_j q_j$ .

These probabilities give the instantaneous transition matrix $\mathbf{Q}$ which we express in the component form. Let $s_{n,q_1,...q_K}$ be the probability of the state $(q_1,\ldots,q_K)$ at the time $n$ . The transition equations of a multiservice loss system have the following form:

$$s_{n+1,q_1,\ldots,q_K}=s_{n,q_1,\ldots,q_K}-\sum_{j=1}^{K}\frac{\lambda_j T}{N}s_{n,q_1,\ldots,q_K}1_{j\in R_{\bar{q}}} \qquad (2)$$

$$-\sum_{j=1}^{K}\frac{\mu_j T}{N}q_j s_{n,q_1,\ldots,q_K}$$

$$+\sum_{j=1}^{K}\frac{\lambda_j T}{N}s_{n,q_1,\ldots,q_j-1,\ldots,q_K}1_{j\in R_{\bar{q}-e_j}}$$

$$\sum_{\bar{q}\in\Omega}s_{n,q_1,...q_K}=1,\ 0\le s_{n,q_1,...,q_K}\le 1.$$

*Remark 1.* For the policies in Examples 1 and 2 the recursion equations (2) satisfy for every $j$ the condition

$$s_{n,q_1,...q_j+1,..q_K}=0 \qquad \text{if } j\notin R_{\bar{q}} . \qquad (3)$$

This is so because if $j\notin R_{\bar{q}}$ , then $\bar{q}+e_j\notin\Omega$ . It might appear that (3) is always true, but even though the transition from $\bar{q}$ to $\bar{q}+e_j$ is blocked by the policy, there could be some other path to get to the state $\bar{q}+e_j$ . Therefore, (3) is a property of the policy.

If the condition in Remark 1 holds, the stationary state solution to the equation set (2) is (see e.g. [8, p. 160] or [9])

$$\pi_q=\pi_{q_1,...q_K}=G^{-1}\prod_{j=1}^{K}\frac{1}{q_j!}\rho_j^{q_j}1_{\bar{q}\in\Omega} \qquad (4)$$

where $G=\sum_{\bar{q}\in\Omega}\prod_{j=1}^{K}\frac{1}{q_j!}\rho_j^{q_j}$ , $\rho_j=\frac{\lambda_j}{\mu_j}$ .

Let us assume that blocking a call of type $j$ creates $\omega_j$ units of cost, and that calls are blocked if there are not enough capacity units. In the stationary state the process (2) generates the average cost $R$ where

$$R=\frac{T}{N}g \ ,\ g=\sum_{\bar{q}\in\Omega}r_{\bar{q}}\pi_{\bar{q}} \ ,\ r_{\bar{q}}=\sum_{j\notin R_{\bar{q}}}\omega_j\lambda_j \ . \qquad (5)$$

The number $g$ is the average cost rate and $r_{\bar{q}}$ is the average cost rate in the state $\bar{q}$ .

## III. Charging by Shadow Prices in a Multiservice Loss System

### A. *Average costs*

The optimal pricing policy in congestion pricing is that the operator charges *shadow prices* – the expectations of future losses from blocked traffic. The shadow prices can be computed in the following way.

The multiservice loss system approaches a steady state and in the steady state it generates constant cost in each time step. The initial state determines how much cost the process generates before reaching the steady state and in most cases rather few steps from the beginning of the process determine quite well how much costs differ when the process is started in a given initial state. This cost difference can be calculated using Howard's equation (see [7])

$$+\sum_{j=1}^{K}\frac{\mu_j T}{N}(q_j+1)s_{n,q_1,\ldots,q_j+}$$

$$\mathbf{QV}=g\mathbf{1}-\mathbf{r} \qquad (6)$$

where $\mathbf{Q}$ is the instantaneous transition matrix of a Markov chain, $\mathbf{V}$ is the vector of relative costs, $g$ is the average cost, and $\mathbf{r}$ is the vector of state cost rates. One method of solving Howard’s equation is inverting $\mathbf{Q}$:

$$\mathbf{V} = \mathbf{Q}^{-1}(g\mathbf{1} - \mathbf{r}). \quad (7)$$

Howard's equation for future relative cost expectations for a multi-service loss system in component form is:

$$\sum_{j \in R_{\bar{q}}} \lambda_j (v(\bar{q}+e_j) - v(\bar{q})) - \sum_{j=1}^{K} \mu_j q_j (v(\bar{q}) - v(\bar{q}-e_j))$$
$$= g - \sum_{j \notin R_{\bar{q}}} \lambda_j \omega_j . \quad (8)$$

The variable $v(\bar{q})$ is called relative cost. The exact solution to Howard's equation is not unique: it is possible to add a constant to $v(\bar{q})$. The state-dependent link *shadow price* for traffic class $k$ is $p_k(\bar{q}) = v_{\bar{q}+e_j} - v_{\bar{q}}$, see [7] formula (7).

If the dimension of the multiservice loss system is small, we can simply invert matrix $\mathbf{Q}$ and solve Howard's equation. If the dimension is very large this method becomes impractical and approximate methods are often proposed. One such method is in [7]. The authors select a small-dimensional basis, compute a projection to the basis, and obtain very good numerical results compared to earlier methods [10], [11], and [12]. We will here denote the capacity $d_j$ in Example 1 by $b_j$. In the special case $b_j = b$, $\mu_j = \mu$ and the policy of Example 1 the exact solution to Howard's equation is (see Appendix I):

$$v(\bar{q}) = \frac{g}{\mu\rho} \sum_{i=1}^{q} \sum_{m=0}^{q-i} \frac{(q-i)!}{(q-i-m)!} \rho^{-m} \quad (9)$$

where $\rho = \sum_{k=1}^{K} \frac{\lambda_k}{\mu_k}$ and $q = \sum_{k=1}^{K} q_k$.

In more general cases we can start with an approximation and complete it into an exact solution in a way described in Appendix I. This way produces a series form solution.

In the case $b_j = b$ but $\mu_j \neq \mu_k$ for some $j$ and $k$, the following closed form approximation seems fairly good

$$\tilde{v}(\bar{q}) = \sum_{j=1}^{K} \frac{q_j}{q} \frac{g}{\mu_j \rho} \sum_{i=1}^{q} \sum_{m=0}^{q-i} \frac{(q-i)!}{(q-i-m)!} \rho^{-m} . (10)$$

In the general case where $b_j \neq b_k$ for some $j$ and $k$ and $\mu_j \neq \mu_k$ for some $j$ and $k$, a rather natural approximation to start from is

$$\hat{v}(\bar{q}) = \sum_{j=1}^{K} \frac{b_j q_j}{c} \left(\frac{b_j}{b}\right)^2 \frac{g}{\mu_j \rho} \sum_{i=1}^{\left[\frac{c}{b_j}\right]} \sum_{m=0}^{\left[\frac{c}{b_j}\right]-i} \frac{\left(\left[\frac{c}{b_j}\right]-i\right)!}{\left(\left[\frac{c}{b_j}\right]-i-m\right)!} \rho_j^{-m}$$

where $[x]$ means the largest integer smaller than $x$

$$b = \sum_{j=1}^{K} b_j , \ \rho_j = \left(\frac{b}{b_j}\right)^2 \rho$$

$$\rho = \sum_{k=1}^{K} \frac{\lambda_k}{\mu_k} \left(\frac{b_k}{b}\right)^2 \text{ and } c = \sum_{k=1}^{K} b_k q_k . \quad (11)$$

These approximations are derived in Appendix I. They are simple to calculate and do not involve any matching of parameters. They are not expected to be as good as in [7] where a large number of free parameters are adjusted to specific numeric values of the loss system. If the shadow prices are used in pricing, we do not need high precision: the input parameters are not known precisely and the price will be rounded to a less precise figure.

In case better accuracy is needed we can calculate more terms in the following way. Any first order approximation $u(\bar{q})$ of $v(\bar{q})$ can be completed to an exact solution of (8) by the following procedure. Let us define recursively

$$f_{1,j}(\bar{q}) = 1 - \Delta_j u(\bar{q}) \quad \text{and for } n \geq 1 , \quad (12)$$

$$f_{n+1,j}(\bar{q}) = -\sum_{\substack{k=1 \\ k \neq j}}^{K} \Delta_k \frac{1}{\mu_j} h(f_{n,j}(\bar{q}), q_j, \rho_j)$$

where $\rho_j = \lambda_j / \mu_j$ and

$$h(f(\bar{q}), q_j, \rho) = \frac{1}{\rho} \sum_{i=1}^{q_j} \sum_{m=0}^{q_j - i} f(\bar{q} - (i+m)e_j) \frac{(q_j - i)!}{(q_j - i - m)!} \rho^{-m} ,$$

$$\Delta = \sum_{j=1}^{K} \Delta_j , \quad (13)$$

$$\Delta_j f(\bar{q}) = \lambda_j (f(\bar{q}+e_j) - f(\bar{q})) - \mu_j q_j (f(\bar{q}) - f(\bar{q}-e_j)) .$$

The function

$$V_{\bar{q}} = g u(q) + g \sum_{j=1}^{K} \sum_{n=1}^{\infty} \frac{1}{\mu_j} h(f_{n,j}(\bar{q}), q_j, \rho_j) \quad (14)$$

satisfies Howard's equation (8). An example of (14) is calculated in Appendix I starting from the approximation (10).

If the operator charges the user with the shadow prices and the network is in the steady state, the distribution of the cost to the user are calculated by multiplying the steady-state probabilities with the shadow prices.

## B. Cost Distribution using Shadow Prices

The relative costs, and consequently the shadow prices, are transient phenomena. After a short time the system will be close to the steady state and producing the average cost rate in (5) regardless of what was the original state. More important parameter than the shadow price is the distribution of the cost in the future. It is not a transient phenomenon and it reflects the differences of using different prices on the services. As the average cost rate depends on the arrival rates $\lambda_i$ and the costs $\omega_i$, we may expect the cost distribution also to depend on these parameters. It is interesting to see if offering more services gives smaller variance in the cost distribution than offering fewer services. We can calculate the cost distribution by solving a set of recursion equations with the initial state of an empty system. If the system is not empty, there comes the transient distribution giving the average costs that are calculated with Howard's equation. These transient costs stay constant and become ignorable in the limit, thus starting from an empty system gives the same limit distribution.

In order to calculate the expected cost distribution up to a finite time horizon $T$ we describe the process using a state vector $(n,\bar{q},r)$ and its probability $s_{n,\bar{q},r}$. Then cost is treated as a state variable and cost distributions can be obtained.

Let $s_{n,\bar{q},r}$ be the probability that the state is $(n,q_1,\ldots,q_K)$ and that the cost vector has the value $r$. Let $r$ be here a nonnegative integer and let us require that each $\omega_j$ is also an integer. Naturally we also require

$$s_{n,q_1,\ldots q_K} = \sum_{r=0}^{\infty} s_{n,q_1,\ldots q_K,r} \ . \qquad (15)$$

We can derive a recursion equation set for the cost distribution by adding the parameter $r$ to (2):

$$s_{n+1,q_1,\ldots,q_K,r} = s_{n,q_1,\ldots,q_K,r} - \sum_{j=1}^{K} \frac{\lambda_j T}{N} s_{n,q_1,\ldots,q_K,r} \qquad (16)$$

$$+ \sum_{j=1}^{K} \frac{\lambda_j T}{N} s_{n,q_1,\ldots,q_K,r-\omega_j} 1_{j \notin R_{\bar{q}}}$$

$$- \sum_{j=1}^{K} \frac{\mu_j T}{N} q_j s_{n,q_1,\ldots,q_K,r}$$

$$+ \sum_{j=1}^{K} \frac{\lambda_j T}{N} s_{n,q_1,\ldots,q_j-1,\ldots,q_K,r} 1_{j \in R_{\bar{q}-e_j}}$$

$$+ \sum_{j=1}^{K} \frac{\mu_j T}{N} (q_j+1) s_{n,q_1,\ldots,q_j+1,\ldots,q_K,r} \ ,$$

$$s_{n,q_1,\ldots q_K} = \sum_{r=0}^{\infty} s_{n,q_1,\ldots q_K,r} \quad , \quad \sum_{\bar{q}\in\Omega} s_{n,q_1,\ldots q_K} = 1 \quad ,$$

$$0 \le s_{n,q_1,\ldots q_K,r} \le 1 .$$

When (16) is summed over $r$, the process (2) is obtained.

The exact solution to (16) is calculated in Appendix II for the policy of Example 2.

*Theorem 1.* Assume that the policy is as in Example 2. Equation set (16) has a limit solution of the form

$$s_{t,\bar{q},r} = \sum_{\substack{r_1=\ldots=r_K=0 \\ r=\sum_{i=1}^{K}\omega_i r_i}}^{\infty} s_{t,r_1,\ldots,r_K} s_{t,\bar{q},r_1,\ldots,r_K} \qquad (17)$$

where

$$s_{t,\bar{q},r_1,\ldots r_K} = \sum_{i=1}^{K} s_{t,\bar{q}_i,r_i} G^{-1} \prod_{\substack{j=1 \\ j\ne i}}^{K} \frac{\rho_j^{q_j}}{q_j!} \qquad (19)$$

$$s_{t,r_1,\ldots r_K} = A \prod_{i=1}^{K} \frac{1}{r_i!} (t\lambda_i)^{r_i} e^{-t\sum_{i=1}^{K}\lambda_i} \qquad (20)$$

$$s_{t,q_i,r_i} = \frac{1}{q_i!} f_i(q_i,r_i,t) g_i(r_i,t) \ ,$$

$G^{-1}$ is from (4). The functions $f_i(q_i,r_i,t)$ and $g_i(r_i,t)$ are complicated and given in Appendix II.

For the policy in Example 1 we may use an approximation

$$C_i = C_i(q) = \max_{q_i} q_i \le \frac{C}{b_i} - \sum_{\substack{m=1 \\ m\ne i}}^{K} \frac{b_m}{b_i} q_m \qquad (21)$$

in Theorem 1.

The formula (20) describes how the cost distribution depends on the arrival rates $\lambda_i$ and service costs $\omega_i$. The solution (17) has an exact form given in Appendix II but it is very complicated. This raises the question whether a simpler pricing scheme could be used instead of shadow prices.

## IV. A Simpler Pricing Scheme

As was seen in Section III, it is difficult to calculate the shadow prices from (8) and the cost distribution from (16). This section presents a simpler pricing scheme that overcomes these difficulties while producing the same average costs.

### *A. Different cost distributions of a multiservice loss system*

One way is to use a product:

$$s_{n,\bar q,r} = s_{n,q_1,\dots q_K,r} = s_{n,r}(\bar q) s_{n,q_1,\dots q_K} \qquad (22)$$

$s_{n,r}(\bar q) = 0$ if $r < 0$,

which satisfies the recursion equation

$$s_{n+1,\bar q,r} = s_{n,\bar q,r} + s_{n,r}(\bar q)\Delta s_{n,\bar q} - \frac{T}{N}\sum_{j\notin R_{\bar q}} \lambda_j s_{n,\bar q,r} + \frac{T}{N}\sum_{j\notin R_{\bar q}} \lambda_j s_{n,\bar q,r-\omega_j}$$

$$s_{n,q_1,\dots q_K} = \sum_{r=0}^{\infty} s_{n,q_1,\dots q_K,r}, \sum_{q\in\Omega} s_{n,q_1,\dots q_K} = 1,$$

$$0 \le s_{n,q_1,\dots q_K,r} \le 1 \qquad (23)$$

The operator $\Delta$ is shorthand to the infinitesimal transition matrix of a multiservice loss system:

$$\Delta f_{n,q_1,\dots q_K} = -\sum_{j\in R_q} \frac{\lambda_j T}{N} f_{n,q_1,\dots,q_K} - \sum_{j=1}^{K} \frac{\mu_j T}{N} q_j f_{n,q_1,\dots,q_K}$$

$$+ \sum_{j\in R_{q-e_j}} \frac{\lambda_j T}{N} f_{n,q_1,\dots,q_j-1,\dots,q_K}$$

$$+ \sum_{j=1}^{K} \frac{\mu_j T}{N}(q_j+1) f_{n,q_1,\dots,q_j+1,\dots,q_K}. \qquad (24)$$

When (23) is summed over $r$, it gives (2). This is seen as follows:

$$s_{n,q_1,\dots q_K} = \sum_{r=0}^{\infty} s_{n,q_1,\dots q_K,r} = s_{n,\bar q} \sum_{r=0}^{\infty} s_{n,r}(\bar q), \qquad (25)$$

thus $\sum_{r=0}^{\infty} s_{n,r}(\bar q) = 1$.

$$\sum_{r=0}^{\infty} s_{n,r-\omega_j}(\bar q) = \sum_{r'=-\omega_j}^{\infty} s_{n,r'}(\bar q) = \sum_{r'=0}^{\infty} s_{n,r'}(\bar q),$$

i.e., $\sum_{r=0}^{\infty}\left(-\sum_{j\notin R_q} \lambda_j s_{n,\bar q,r} + \sum_{j\notin R_q} \lambda_j s_{n,\bar q,r-\omega_j}\right) = 0$ and we obtain (2) in the form:

$$s_{n+1,\bar q} = \sum_{r=0}^{\infty} s_{n+1,\bar q,r} = s_{n,\bar q} + \Delta s_{n,\bar q}. \qquad (26)$$

Another way is to take a slightly more general form:

$$s_{n,\bar q,r} = a s_{n,r}(\bar q) s_{n,q_1,\dots q_K} + (1-a) v_{n,q_1,\dots q_K,r} \qquad (27)$$

where $a \in [0,1]$ and

$$\lim_{n\to\infty}(v_{n+1,\bar q,r} - v_{n,\bar q,r}) = 0 \quad , \quad \sum_{r=0}^{\infty} v_{n,\bar q,r} = s_{n,\bar q} \quad ,$$

$$\lim_{n\to\infty} \frac{v_{n,\bar q,r}}{s_{n,r}(\bar q)} = 0.$$

The corresponding recursion equation

$$s_{n+1,\bar q,r} = s_{n,\bar q,r} + a s_{n,r}(\bar q)\Delta s_{n,\bar q} + (1-a)\Delta v_{n,\bar q,r} \qquad (28)$$

$$-\frac{T}{N}\sum_{j\notin R_{\bar q}} \lambda_j s_{n,\bar q,r} + \frac{T}{N}\sum_{j\notin R_{\bar q}} \lambda_j s_{n,\bar q,r-\omega_j}$$

$$s_{n,q_1,\dots q_K} = \sum_{r=0}^{\infty} s_{n,q_1,\dots q_K,r} \quad , \quad \sum_{q\in\Omega} s_{n,q_1,\dots q_K} = 1 \quad ,$$

$$0 \le s_{n,q_1,\dots,q_K,r} \le 1.$$

also gives (2) when summed over $r$. To see this, let us notice that the last two terms of (28) cancel when summed over $r$ and that

$$\sum_{r=0}^{\infty} \Delta v_{n,\bar q,r} = \Delta \sum_{r=0}^{\infty} v_{n,\bar q,r} = \Delta s_{n,\bar q} = 0. \qquad (29)$$

Then $a\sum_{r=0}^{\infty} s_{n,r}(\bar{q})=1$ because the term

$\sum_{r=0}^{\infty}(1-a)v_{n,\bar{q},r}$ disappears in the limit $n\to\infty$.

*Lemma 1.* The set of equations (16), (23) and (28) all generate the cost rate $R$ in (5).

The proposed simpler pricing scheme is pricing the connections with the future expectations of the process (23) instead of by the future expectations of the process (16). In shadow prices (16) all states $\bar{q}$ are associated with a positive shadow price because there is some future expectation that the system starting from $\bar{q}$ will be blocking traffic in the future, and consequently future expectations of costs are positive. In (23) only such states that block give any charge. We will show this after Theorems 2 and 3 below. While shadow prices are optimal in a certain sense, the differences between shadow prices and pricing by future expectations of (23) are small. User's bills are in any case rounded, user traffic is measured, network congestion state is estimated and in reality shadow prices cannot be very exact.

## B. Cost distribution for the simpler pricing scheme

The recursion equations (23) can be solved exactly because they are created to satisfy the product form solution (22). The limit distribution to the recursion equations (28) can be solved similarly. The part $v_{n,\bar{q},r}$ disappears in the limit.

*Theorem 2.* Assume that the process (23) satisfies $s_{0,\bar{q},r}=0$ if $r\neq 0$. The solution to (23) has the form

$$s_{n,q_1,\dots q_K,r}=\left(\sum_{\substack{r_1=\dots=r_K=0\\(r_1,\dots r_K)\in A_{r,\bar{q}}}}^{n} C_{n,q_1,\dots q_K,r}\right)\pi_{n,\bar{q}}$$

(30)

where

$$C_{n,q_1,\dots q_K,r}=\begin{pmatrix} & n & & \\ r_1 & . & . & . & r_K\end{pmatrix}\left(1-\frac{T}{N}\sum_{i\notin R_{\bar{q}}}\lambda_i\right)^{n-\sum_{i=1}^{K}r_i}\prod_{i\notin R_{\bar{q}}}\left(\frac{T}{N}\lambda_i\right)^{r_i}$$

$$A_{r,\bar{q}}=\left\{(r_1,\dots,r_K)\left|\sum_{i=1}^{K}\omega_i r_i=r,\quad r_i=0\quad\text{if}\quad i\in R_{\bar{q}}\right.\right\}$$

and $\pi_{n,\bar{q}}\to\pi_{\bar{q}}$ when $n\to\infty$, $\pi_{\bar{q}}$ is given by (4).

It is also possible to derive the continuous limit:

*Theorem 3.* The continuous time limit function for (30) when $N\to\infty$, $n=xN$, $t=xT$ is

$$s_{t,q_1,\dots q_K,r}=\left(\sum_{\substack{r_1=\dots=r_K=0\\(r_1,\dots r_K)\in A_{r,\bar{q}}}}^{\infty}\prod_{i\notin R_{\bar{q}}}\frac{1}{r_i!}(t\lambda_i)^{r_i}\right)e^{-t\sum_{i\notin R_{\bar{q}}}\lambda_i}\pi_{\bar{q}} \quad .$$

(31)

Let us note that (30) is not a solution to (16). The formula (30) has the form of a product

$$s_{n,q_1,\dots.q_K,r}=s_{n,r}(\bar{q})\pi_q \,. \qquad (32)$$

Then the detailed balance conditions for each $j$ and $r$ in the limit $n\to\infty$

$$\lim_{n\to\infty}\mu_i(q_i+1)s_{n,q_1,\dots q_j+1,\dots q_K,r}=\lim_{n\to\infty}\lambda_i s_{n,q_1,\dots q_K,r} \qquad (33)$$

and

$$\lim_{n\to\infty}\mu_i q_i s_{n,q_1,\dots q_K,r}=\lim_{n\to\infty}\lambda_i s_{n,q_1,\dots q_i-1,\dots q_K,r}$$

do not hold in (16), since $\pi_q$ satisfies them, but if for some $q$ holds $R_q\neq R_{q+e_i}$, then $s_{n,r}(q+e_i)\neq s_{n,r}(q)$.

Let us check that (31) produces the correct cost.

*Lemma 2.* The solution (31) satisfies

$$\sum_{r=0}^{\infty} r s_{t,\bar{q},r}=tr_{\bar{q}}\pi_{\bar{q}} \qquad (34)$$

*Remark 2.* Let us calculate directly two easy examples:

a) If $R_q = \{1,2\}$ and $\omega_1 = \omega_2 = 1$, then

$$\sum_{r=0}^{\infty} rs_{t,q_1,\ldots q_K,r}$$

$$= \sum_{r=0}^{\infty} r \sum_{r_1=0}^{r} \frac{1}{r_1!}\frac{1}{(r-r_1)!}(xT\lambda_1)^{r_1}(xT\lambda_2)^{r-r_1} e^{-xT(\lambda_1+\lambda_2)}\pi_{\bar{q}}$$

$$= \sum_{r=0}^{\infty} \frac{1}{(r-1)!}(xT)^r(\lambda_1+\lambda_2)^r e^{-xT(\lambda_1+\lambda_2)}\pi_{\bar{q}}$$

$$= xT(\lambda_1+\lambda_2)\pi_{\bar{q}} = tr_{\bar{q}}\pi_{\bar{q}}. \qquad (35)$$

b) If $R_q = \{1\}$, then the multinomial recurrence formula

$$\begin{pmatrix} & n-1 & & \\ i_1 & . & . & . & i_m \end{pmatrix} = \left(1 - \frac{1}{n}\sum_{k=1}^{m} i_k\right)\begin{pmatrix} & n & & \\ i_1 & . & . & . & i_m \end{pmatrix}$$

for the binomial term gives for $n >> 1$:

$$\sum_{r=0}^{\infty} rs_{n,q_1,\ldots q_K,r} = \sum_{i_i=0}^{n} \omega_i i_i \binom{n}{i_i}\left(1-\frac{T}{N}\lambda_1\right)^{n-i_1}\left(\frac{T}{N}\lambda_1\right)^{i_1}\pi_{\bar{q}}$$

$$= n\frac{T}{N}\omega_1\lambda_1\pi_{\bar{q}} \qquad (36)$$

Let us verfy that in Theorems 2 and 3 only such states that block give any charge: If the system is in a state that does not block, there is only the point $(r_1,\ldots,r_K) = (0,\ldots,0)$ in $A_{r,\bar{q}}$. Thus, $s_{n,q,r} = 0$ if $r \neq 0$.

Let us notice that the part corresponding to the arrival rates $\lambda_i$ and service costs $\omega_i$ is similar in (17) and (31), only in (17) there is the scaling factor $A$.

## V. CONCLUSIONS

The user and the operator are interested in the variability of the cost of a connection caused by congestion pricing. The distribution of the cost of a connection to a user can be calculated in the ways given in Section III.

An operator should know the cost distribution because it may have to add a price component to account for the risk of the actual costs exceeding the average costs. An exact solution for the distribution of costs to an operator over a longer time period can be obtained for a particular way of distributing costs to users, described by the recursion equation set (23). It is a slightly different pricing strategy than using shadow prices, but it gives the same incentives for users to lower their traffic demands. For shadow prices an exact expression could only be obtained for the simple blocking policy in Example 2. However, the structure of the cost distribution is very similar in (17) and (32), and (32) can be used as an approximation.

The pricing scheme in Theorems 2 and 3 does not much differ from the shadow pricing scheme. Pricing schemes that approximate shadow prices but are more simple have earlier been suggested by other authors, for instance [6] proposes a static pricing scheme. Exact shadow prices are found too complicated and simpler approximations are acceptable substitutes.

Theorems 1, 2 and 3 use the product form solution (4), which gives an exact solution to blocking probabilities for a certain class of routing policies in a circuit-switched network. The normalization constant $G$ is difficult to calculate but effective approximations do exist, or instance [14] and [15]. The expressions in Theorems 2 and 3 are rather similar to that of $G$ and similar effective calculation methods apply. As the cost distribution is only used for risk estimation and not for satisfying Quality of Service requirements, less precise methods are sufficient. Calculation of $G$ precisely is a NP-complete problem and similar complexity is expected to the $r$ dependent part in (17) and (31), but as only the general shape of the distribution is needed, Monte Carlo simulation can give sufficient approximations for practical purposes in a reasonable time.

Theorems 1, 2 and 3 deal with the risk the operator is taking with congestion pricing. There is a general feeling that some price component to cover the risks is needed, but the case of risk in congestion pricing has not been much analyzed in the literature. The general aspect of risk is much wider than is treated in this paper. Many issues, such as growing demand and various aspects of uncertainty [16], [17] and [18] must be omitted in this paper.

## APPENDIX I: PROOFS OF SECTION II A.

*Proof of (9):* Let $v(\bar{q})$ be as in (9), then

$$v(\bar{q}) - v(\bar{q}-e_j) = \frac{g}{\mu\rho}\sum_{m=0}^{q-1}\frac{(q-1)!}{(q-1-m)!}\rho^{-m}. \qquad (37)$$

In the policy of Example 1 holds $1_{c<C} = 1_{j\in R_{\bar{q}}}$, thus

$$\sum_{j=1}^{K} 1_{c<C}\lambda_j(v(\bar{q}+e_j)-v(\bar{q})) - \sum_{j=1}^{K}\mu q_j(v(\bar{q})-v(\bar{q}-e_j))$$

$$= \frac{g}{\mu\rho}\sum_{i=1}^{K}1_{c<C}\lambda_i\sum_{m=0}^{q}\frac{q!}{(q-m)!}\rho^{-m} - \sum_{i=1}^{K}q_i\frac{g}{\rho}\rho\sum_{m=1}^{q}\frac{(q-1)!}{(q-m)!}\rho^{-m}$$

$$= \begin{cases} g & \text{if} \quad c < C \\ A & \text{if} \quad c = C \end{cases}$$

(38)

where

$$A = -g\sum_{m=1}^{C/b}\frac{(C/b)!}{(C/b-m)!}\rho^{-m} \quad (39)$$

There is only one nonzero value for $\sum_{j\notin R_{\bar q}}\lambda_j\omega_j$ in this case and the value of $g$ in (5) sets $A$ to the value satisfying (8):

$$A = g - \sum_{j\notin R_{\bar q}}\lambda_j\omega_j \;. \;\square \quad (40)$$

*Motivation of (10) and (11):* We show why the heuristic approximations (10) and (11) are natural and can give rather good results. Let us assume that at the axis $\bar q = qe_j$ the shadow price has the approximation

$$p_k(\bar q) = v(\bar q + e_k) - v(\bar q) \approx a_k\frac{g}{\mu\rho}\sum_{m=0}^{q_j}\frac{q_j!}{(q_j-m)!}\rho^{-m} \quad . \quad (41)$$

This is a generalization of the relative cost given by (9) and $a_k = 1$ if every $\mu_j = \mu$. In the area when every call is accepted the recursion (8) gives for $\bar q = qe_j$

$$\sum_{k=1}^{K}\lambda_k p_k(\bar q) - \mu_j q_j p_j(\bar q - e_j) = g \quad . \quad (42)$$

For the approximation (11) we get

$$\sum_{k=1}^{K}\lambda_k a_k\frac{g}{\mu\rho}\sum_{m=0}^{q_j}\frac{q_j!}{(q_j-m)!}\rho^{-m} - \mu_j q_j a_j\frac{g}{\mu\rho}\sum_{m=0}^{q_j-1}\frac{(q_j-1)!}{(q_j-1-m)!}\rho^{-m}$$

$$= \sum_{k=1}^{K}\lambda_k a_k\frac{g}{\mu\rho} + \left(\sum_{k=1}^{K}\lambda_k a_k\frac{g}{\mu\rho} - \mu_j a_j\frac{g}{\mu}\right)\sum_{m=0}^{q_j}\frac{q_j!}{(q_j-m)!}\rho^{-m} = g$$

.

This is satisfied if

$$\rho = \sum_{k=1}^{K}\lambda_k\frac{a_k}{\mu} \quad \text{and} \quad 1 = \mu_j a_j\frac{1}{\mu} \;,$$

i.e., $a_j = \dfrac{\mu}{\mu_j}$ , $\rho = \sum_{k=1}^{K}\dfrac{\lambda_k}{\mu_k}$ where $\mu = \sum_{j=1}^{K}\mu_j$ . (43)

*Proof of (14):*

The following lemma holds.

*Lemma 3:*

$$\Delta_j\frac{1}{\mu_j}h(f(\bar q), q_j, \rho_j) = f(\bar q) \quad \text{for } \rho_j = \lambda_j/\mu_j \;.$$

*Proof:*

$$\rho_j(h(f(\bar q + e_j), q_j + 1, \rho_j) - h(f(\bar q), q_j, \rho_j))$$
$$- q_j(h(f(\bar q), q_j, \rho_j) - h(f(\bar q - e_j), q_j - 1, \rho_j))$$

$$= \sum_{i=0}^{q_j}\sum_{m=0}^{q_j-i} f(\bar q - (i+m)e_j)\frac{(q_j-i)!}{(q_j-i-m)!}\rho_j^{-m}$$
$$- \sum_{i=1}^{q_j}\sum_{m=0}^{q_j-i} f(\bar q - (i+m)e_j)\frac{(q_j-i)!}{(q_j-i-m)!}\rho_j^{-m}$$

$$- q_j\frac{1}{\rho_j}\left(\sum_{i=1}^{q_j}\sum_{m=0}^{q_j-i} f(\bar q - (i+m)e_j)\frac{(q_j-i)!}{(q_j-i-m)!}\rho_j^{-m}\right)$$
$$+ q_j\frac{1}{\rho_j}\left(\sum_{i=2}^{q_j}\sum_{m=0}^{q_j-i} f(\bar q - (i+m)e_j)\frac{(q_j-i)!}{(q_j-i-m)!}\rho_j^{-m}\right)$$

$$= \sum_{m=0}^{q_j} f(\bar q - me_j)\frac{(q_j)!}{(q_j-m)!}\rho_j^{-m}$$
$$- q_j\frac{1}{\rho_j}\sum_{m=0}^{q_j-1} f(\bar q - (m+1)e_j)\frac{(q_j-1)!}{(q_j-1-m)!}\rho_j^{-m}$$

$$= \sum_{m=0}^{q_j} f(\bar q - me_j)\frac{(q_j)!}{(q_j-m)!}\rho_j^{-m}$$
$$- \sum_{m=1}^{q_j} f(\bar q - me_j)\frac{(q_j)!}{(q_j-m)!}\rho_j^{-m} = f(\bar q) \;. \;\square \quad (44)$$

By Lemma 3, $V_{\bar q}$ satisfies

$$\Delta V_{\bar q} = g\sum_{j=1}^{K}\Delta_j u(q) + g\sum_{j=1}^{K}\sum_{n=1}^{\infty}\sum_{k=1}^{K}\Delta_k\frac{1}{\mu_j}h(f_{n,j}(\bar q), q_j, \rho_j)$$

$$= g\sum_{j=1}^{K}\left(1 - f_{1,j}(\bar q)\right) + g\sum_{j=1}^{K}\sum_{n=1}^{\infty}\left(f_{n,j}(\bar q) - f_{n+1,j}(\bar q)\right)$$

(45)

$$= g\sum_{j=1}^{K}\left(1 - f_{1,j}(\bar q)\right) + g\sum_{j=1}^{K} f_{1,j}(\bar q) = g \;.$$

The value of $v(\overline{q})$ at the points when some class is blocked is already determined by the values of $v(\overline{q})$ at the points when no class is blocked. This means that we do not need to modify $V_{\overline{q}}$ so that it satisfies (8). Changing the operator $\Delta$ to

$$\Delta = \sum_{j=1}^{K} \Delta_j \, , \qquad (46)$$

$$\Delta_j f(\overline{q}) = 1_{j \in R_{\overline{q}}} \lambda_j (f(\overline{q} + e_j) - f(\overline{q})) - \mu_j q_j (f(\overline{q}) - f(\overline{q} - e_j))$$

the solution $V_{\overline{q}}$ necessarily satisfies (8):

$$\Delta V_{\overline{q}} = g - \sum_{j \notin R_{\overline{q}}} \lambda_j \omega_j \, . \qquad (47)$$

This is so because before states block, the recursion equation is the same as in (13) and therefore the solution agrees with $V_{\overline{q}}$ for those values. At the border when states block, $V_{\overline{q}}$ is already determined by the values at unblocking states. Therefore Howard's equation must be satisfied at the blocking states.

*As an example of (14)* let us complete the approximation $\tilde{v}(\overline{q})$ in (10) to an exact solution of (8). Let us select

$$u(\overline{q}) = \frac{1}{\mu} h(1, q, \rho) \qquad (48)$$

where $q = \sum_{j=1}^{K} q_j$, $\mu = \sum_{j=1}^{K} \mu_j$, $\rho = \sum_{j=1}^{K} \frac{\lambda_j}{\mu_j}$ and

$$h(f(q), q, \rho) = \frac{1}{\rho} \sum_{i=1}^{q} \sum_{m=0}^{q-i} f(q-i-m) \frac{(q-i)!}{(q-i-m)!} \rho^{-m} \, .$$

Let $\Delta$ be as in (46). A similar calculation as in the proof of Lemma 3 shows that

$$\Delta \frac{1}{\mu} h(1, q, \rho) = \sum_{j=1}^{K} \frac{\lambda_j}{\mu} \frac{1}{\rho} \sum_{m=0}^{q} \frac{q!}{(q-m)!} \rho^{-m} - \sum_{j=1}^{K} \frac{\mu_j}{\mu} q_j \sum_{m=1}^{q} \frac{(q-1)!}{(q-m)!} \rho^{-m} \, . \qquad (49)$$

Then

$$f_{1,j}(\overline{q}) = \left( \frac{\lambda_j}{\mu_j} - \frac{\lambda_j}{\mu} \right) \frac{1}{\rho} \sum_{m=0}^{q} \frac{q!}{(q-m)!} \rho^{-m} - \left( 1 - \frac{\mu_j}{\mu} \right) q_j \sum_{m=1}^{q} \frac{(q-1)!}{(q-m)!} \rho^{-m} \, . \qquad (50)$$

At any coordinate axis $\overline{q} = q e_j$ holds $\tilde{v}(q) = g u(\overline{q})$. The complete solution to (8) is

$$V_{\overline{q}} = g \frac{1}{\mu} h(1, q, \rho) + g \sum_{j=1}^{K} \sum_{n=1}^{\infty} \frac{1}{\mu_j} h(f_{n,j}(\overline{q}), q_j, \rho_j) \, .$$

## APPENDIX II: PROOF OF THEOREM 1

For notational simplicity, we write $\overline{q}$ simply as $q$. We only intend to prove that (16) has a solution of the form (17) and do this by constructing such a solution. The solution turns out to be very complicated and not of practical usage. However, Lemma 4 below is interesting in its own right. The equation is reoccurring often in Markov chains with a cost parameter and the solution is not that easy to derive.

*Lemma 4.* Let us assume that $f : \mathbf{N} \to \mathbf{R}$ satisfies the recursion equation

$$f(q+2) - (q+a) f(q+1) + \rho(q+1) f(q) = 0 \, , \quad f(q) = 0 \text{ if } q < 0 \, , \qquad (51)$$

where $q \in \mathbf{Z}$, $\rho, a \in \mathbf{R}$, $a - \rho \notin \mathbf{N}$ and $\rho \neq 0$.

Then the function $f$ is

$$f(q) = \sum_{j=1}^{q + \lfloor a-\rho \rfloor} \gamma_j \Gamma(q + a - \rho - j) + C \sum_{j=1}^{q + \lceil a-\rho \rceil} \gamma_j \Gamma(q + a - \rho - j) \qquad (52)$$

where

$$\alpha_1 = \gamma_1 \rho(\rho + 2 - a) \, ,$$

$\alpha_{i+1} = \frac{\rho}{i} \alpha_i (2 + \rho - a + i)$, $i \geq 1$, and $\gamma_{i+1} = \frac{\alpha_i}{i}$, $i \geq 1$ and $C$ is a real number.

*Proof*: Let us define

$$f_i(q) = \sum_{j=1}^{i} \gamma_j \Gamma(q + a - \rho - j)$$

where $\Gamma(x+1) = \prod_{j=1}^{\lfloor x \rfloor} (j + x - \lfloor x \rfloor)$. (53)

Write

$$A_i = f_i(q+2) - (q+a) f_i(q+1) + \rho(q+1) f_i(q)$$

.

Then

$$A_{i+1} - A_i = \gamma_{i+1}\Gamma\ (q+2+a-\rho-i-1) \quad (54)$$

$$-(q+a)\gamma_{i+1}\Gamma(q+1+a-\rho-i-1)$$

$$+\rho(q+1)\gamma_{i+1}\Gamma(q+a-\rho-i-1).$$

From $f_1(q) = \gamma_1 \Gamma(q+a-\rho-1)$ we get the initial step to the induction:

$$A_1 = \gamma_1 \rho(\rho+2-a) f_1(q) = \alpha_1 \Gamma(q+a-\rho-1)$$

when we define $\alpha_1 = \gamma_1 \rho(\rho+2-a)$. Let the induction assumption be that for $j = 1,...,i$ and for some numbers $\alpha_j$ holds

$$A_j = \alpha_j \Gamma(q+a-\rho-j) \quad (55)$$

Then for $j = i+1$ we get

$$A_{i+1} - A_i$$

$$= \gamma_{i+1}\Gamma(q+a-\rho-i-1)((q+a-\rho-i)(q+a-\rho-i-1)$$

$$-(q+a)(q+a-\rho-i-1) + \rho q + \rho) \quad (56)$$

$$= \gamma_{i+1}\Gamma(q+a-\rho-i-1)(-iq-(a-\rho-i-1)(\rho+i)+\rho)$$

.

By the induction assumption

$$A_i = \alpha_i (q+a-\rho-i-1)\Gamma(q+a-\rho-i-1).$$

Then

$$A_{i+1} = \Gamma(q+a-\rho-i-1)((\alpha_i - \gamma_{i+1} i)q$$

$$+\alpha_i(a-\rho-i-1) - \gamma_{i+1}(a-\rho-i-1)(\rho+i) + \gamma_{i+1}\rho)$$

$$= \alpha_{i+1}\Gamma(q+a-\rho-i-1) \quad (57)$$

when we choose

$$\gamma_{i+1} = \frac{\alpha_i}{i} \quad \text{and} \quad \alpha_{i+1} = \frac{\rho}{i}\alpha_i(2+\rho-a+i). \quad (58)$$

By induction we can continue so long that $q$ is removed from the Gamma function and we have a scalar. The scalar is unfortunately not usually the correct initial value. Let us define

$$A_{\lfloor a-\rho \rfloor + q} = \alpha_{\lfloor a-\rho \rfloor + q}\Gamma(a-\rho-\lfloor a-\rho \rfloor) \quad \text{and}$$

$$A_{\lceil a-\rho \rceil + q} = \alpha_{[a-\rho]+q}\Gamma(a-\rho-[a-\rho]). \quad (59)$$

These are scalars, i.e., not functions of $q$. If $a-\rho \notin \mathbf{N}$, we get two different scalars that can be weighted to that their linear combination is the desired initial value (zero). Thus

(60)

$$f(q) = \sum_{j=1}^{q+\lfloor a-\rho \rfloor} \gamma_j \Gamma(q+a-\rho-j) + C \sum_{j=1}^{q+\lceil a-\rho \rceil} \gamma_j \Gamma(q+a-\rho-j)$$

satisfies (51) if $C = -A_{\lfloor q+a-\rho \rfloor} / A_{\lceil q+a-\rho \rceil}$.The parameter $\gamma_1$ is a free parameter and we can set $\gamma_1 = 1$. □

*Proof of Theorem 1:* Let us use the trial

$$s_{n,q,r} = \sum_{\substack{r_1=...=r_K=0 \\ r=\sum_{i=1}^{K}\omega_i r_i}}^{\infty} s_{n,q,r_1,...,r_K} \quad (61)$$

and require that each term of the sum satisfies (16). Using the property (3), we can change the position of $1_{j \in R_{\bar{q}}}$ and rewrite (16) in the following form:

$$s_{n+1,q,r_1,\ldots,r_K} = s_{n,q,r_1,\ldots,r_K} \quad (62)$$

$$-\frac{T}{N}\sum_{i=1}^{K}\mu_i q_i s_{n,q,r_1,\ldots,r_K} - \frac{T}{N}\sum_{i=1}^{K}\lambda_i s_{n,q,r_1,\ldots,r_K}$$

$$+\frac{T}{N}\sum_{i=1}^{K}\mu_i(q_i+1)s_{n,q+e_i,r_1,\ldots,r_K} 1_{i \in R_q}$$

$$+\frac{T}{N}\sum_{i=1}^{K}\lambda_i s_{n,q-e_i,r_1,\ldots,r_K}$$

$$+\frac{T}{N}\sum_{i=1}^{K}\lambda_i s_{n,q,r_1,\ldots r_i-1,\ldots,r_K} 1_{i\notin R_q} .$$

Let us look for a solution of the type

$$s_{n,r_1,\ldots r_K} s_{t,q,r_1,\ldots,r_K} = s_{n,r_1,\ldots r_K} \sum_{i=1}^{K} s_{t,q_i,r_i} G^{-1} \prod_{\substack{j=1\\ j\neq i}}^{K} \frac{\rho_j^{q_j}}{q_j!}$$

(63)

where

$$s_{n+1,r_1,\ldots,r_K} = s_{n,r_1,\ldots,r_K} - \frac{T}{N}\sum_{i=1}^{K} a_i\lambda_i s_{n,r_1,\ldots,r_K} + \frac{T}{N}\sum_{i=1}^{K} a_i\lambda_i s_{n,r_1,\ldots,r_i-1,\ldots r_K}$$

$$s_{0,r_1,\ldots,r_K} = \delta_{r_1=\ldots=r_K=0} . \qquad (64)$$

Summing over $r_i$ gives

$$s_{n,,r} = \left( \sum_{\substack{r_1=\ldots=r_K=0\\ r=\sum_{i=1}^{K}\omega_i r_i}}^{\infty} \binom{n}{r_1 \;\; . \;\; . \;\; . \;\; r_K} \left(1-\frac{T}{N}\sum_{i=1}^{K} a_i\lambda_i\right)^{n-\sum_{i=1}^{K} r_i} \prod_{i=1}^{K}\left(\frac{T}{N}a_i\lambda_i\right)^{r_i} \right)$$

.

Taking the continuous limit of this expression gives (20). Let us modify (62) by adding and removing terms:

$$s_{n+1,r_1,\ldots,r_K} s_{q,r_1,\ldots r_K} = s_{n,r_1,\ldots,r_K} s_{q,r_1,\ldots,r_K} \qquad (65)$$

$$-\frac{T}{N}\sum_{i=1}^{K} a_i\lambda_i s_{n,r_1,\ldots,r_K} s_{q,r_1,\ldots r_K} + \frac{T}{N}\sum_{i=1}^{K} a_i\lambda_i s_{n,r_1,\ldots,r_i-1,\ldots r_K} s_{q,r_1,\ldots r_K}$$
$$+\frac{T}{N}\sum_{i=1}^{K} a_i\lambda_i s_{n,r_1,\ldots,r_K} s_{q,r_1,\ldots r_K} - \frac{T}{N}\sum_{i=1}^{K} a_i\lambda_i s_{n,r_1,\ldots,r_i-1,\ldots r_K} s_{q,r_1,\ldots r_K}$$
$$-\frac{T}{N}\sum_{i=1}^{K}\mu_i q_i s_{n,q,r_1,\ldots,r_K} - \frac{T}{N}\sum_{i=1}^{K}\lambda_i s_{n,q,r_1,\ldots,r_K}$$

$$+\frac{T}{N}\sum_{i=1}^{K}\mu_i(q_i+1) s_{n,q+e_i,r_1,\ldots,r_K} 1_{i\in R_q} + \frac{T}{N}\sum_{i=1}^{K}\lambda_i s_{n,q-e_i,r_1,\ldots,r_K}$$

$$+\frac{T}{N}\sum_{i=1}^{K}\lambda_i s_{n,q,r_1,\ldots r_i-1,\ldots,r_K} 1_{i\notin R_q} .$$

The first four terms give zero because of (64), thus the other terms also sum to zero. Let us take the continuous limit $n/N \to t/T$ and define

$$h_i(r_i,t) = \lim_{n\to t} \frac{s_{n,r_1,\ldots r_i-1,\ldots r_K}}{s_{n,r_1,\ldots r_K}} = \frac{r_i}{ta_i\lambda_i} , \qquad (66)$$

where $t$ is the continuous time parameter. Dividing (65) by $s_{n,r_1,\ldots r_K}$ and inserting $h_i(r_i,t)$ and

$$s_{t,q_i,r_i} = \frac{1}{q_i!} f_i(q_i,r_i,t) g_i(r_i,t)$$

(67)

gives

$$0 = \sum_{i=1}^{K}\prod_{\substack{j=1\\ j\neq i}}^{K} \frac{\rho_j^{q_j}}{q_j!}\lambda_i\theta_i$$

where

$$\theta_i = (a_i - 1 - a_i h_i(r_i,t))\frac{1}{q_i!} f_i(q_i,r_i,t) g_i(r_i,t)$$
$$-\rho_i^{-1}\frac{1}{(q_i-1)!} f_i(q_i,r_i,t) g_i(r_i,t)$$

(68)

$$+\rho_i^{-1}\frac{1}{q_i!} f_i(q_i+1,r_i,t) g_i(r_i,t) 1_{i\in R_q}$$

$$+\frac{1}{(q_i-1)!} f_i(q_i-1,r_i,t) g_i(r_i,t)$$

$$+h_i(r_i,t)\frac{1}{q_i!} f_i(q_i,r_i-1,t) g_i(r_i-1,t) 1_{i\notin R_q} .$$

Let us satisfy (68) by setting the $\theta_i$ to zero for every $i$, thus

$$0 = f_i(q_i+1,r_i,t) g_i(r_i,t) 1_{i\in R_q} \qquad (69)$$

$$-(q_i+\rho_i(1-a_i+a_i h_i(r_i,t))) f_i(q_i,r_i,t) g_i(r_i,t)$$

$$+\rho_i q_i f_i(q_i-1,r_i,t) g_i(r_i,t)$$

$$+\rho_i h_i(r_i,t) f_i(q_i,r_i-1,t) g_i(r_i-1,t) 1_{i\notin R_q} .$$

If $i \in R_q$ we get

$$f_i(q_i+1,r_i,t) - (q_i+\rho_i(1-a_i+a_i h_i(r_i,t))) f_i(q_i,r_i,t)$$
$$+\rho_i q_i f_i(q_i-1,r_i,t) = 0 . \qquad (70)$$

This is a recursion only for $q_i$, while $r_i$ and $t$ are only parameters. Lemma 4 gives the solution to (70). For $i \notin R_q$ in the policy of Example 2 we can satisfy (70) by assigning $q_i = C_i$. Here $h_i(r_i, t) \neq 0$ and thus (69) gives a recursion equation from which $g_i(r_i)$ is solved by the form

$$g_i(r_i, t) = g_i(r_i - 1, t) \prod_{j=1}^{r_i} D_{i,r_i} ,$$

$$D_{i,j} = \frac{\rho_i h_i(j,t) f_i(C_i, j-1, t)}{C_i + \rho_i(1 - a_i + a_i h_i(j,t)) f_i(C_i, j, t) - \rho_i f_i(C_i - 1, j, t)}$$

where $C_i$ is the maximum value for $q_i$.

Let us look at the special case $h_i(r_i, t) = 0$. Then $f_i(q_i, r_i, t) = f_i(r_i)$. Instead of (67) we will select

$$s_{t,q_i,r_i} = \frac{1}{q_i!} f_i(q_i) g_i(r_i) t^{-r_i} \quad (72)$$

and in (68)-(69)

$$h_i(r_i) \frac{t^{-(r_i-1)}}{t^{-r_i}} = \frac{r_i}{a_i \lambda_i} \neq 0 \quad . \quad (73)$$

We get an equation similar to (69) and can solve $g_i(r_i)$. Scaling by selecting the parameter $A$ lets us satisfy

$$\sum_{r,q} s_{n,q,r} = 1 \text{ and } 1 \geq s_{n,q,r} \geq 0 . \quad (74)$$

There is still the requirement that the sum over $r$ of values $s_{n,q,r}$ must be finite. The part $s_{n,r}$ in $s_{n,q,r}$ satisfies it. □

## APPENDIX III: PROOFS OF SECTION IV

For notational simplicity, we write $\bar{q}$ simply as $q$.

*Proof of Lemma 1:* Since (23) is a special case of (28), it suffices to show the claim for (28) and (16). In (28) holds

$$\lim_{n\to\infty} \sum_{r=0}^{\infty} r(as_{n,r}(q)\Delta s_{n,q} + (1-a)\Delta v_{n,q,r}) = 0$$

because $\lim_{n\to\infty} \Delta s_{n,q} = \Delta \pi_q = 0$ and $\lim_{n\to\infty} \Delta v_{n,q,r} = 0$.

Thus

$$\sum_{r=0}^{\infty} r\left(s_{n+1,q,r} - s_{n,q,r}\right) = \frac{T}{N} \sum_{j\notin R_q} \lambda_j \omega_j . \quad (75)$$

To prove the claim to (16), it is only needed to multiply (16) by $r$ and sum over $r$

$$\sum_{r=0}^{\infty} r\left(s_{n+1,q,r} - s_{n,q,r}\right) = \frac{T}{N} \sum_{j\notin R_q} \lambda_j \omega_j . \square \quad (76)$$

*Proof of Theorem 2:* It is sufficient to show that the equations

$$\pi_{q_1,\ldots q_K} = \lim_{n\to\infty} \sum_{r=0}^{\infty} s_{n,q_1,\ldots q_K,r} \quad , \quad 0 \leq s_{n,q_1,\ldots q_K,r} \leq 1 \quad , \quad (77)$$

$$s_{n+1,q_1,\ldots q_K,r} = s_{n,q_1,\ldots q_K,r} - \frac{T}{N} \sum_{j\notin R_q} \lambda_j s_{n,q_1,\ldots,q_K,r} + \frac{T}{N} \sum_{j\notin R_q} \lambda_j s_{n,q_1,\ldots,q_K,r-\omega_j}$$

are satisfied. Let us write

$$s_{n,r_1,\ldots r_K}(q) = \binom{n}{r_1 \; . \; . \; . \; r_K} \left(1 - \frac{T}{N} \sum_{i\notin R_q} \lambda_i\right)^{n - \sum_{i=1}^{K} r_i} \prod_{i\notin R_q} \left(\frac{T}{N}\lambda_i\right)^{r_i} \quad (78)$$

and re-index

$$\left\{ j \middle| j \notin R_q \right\} = \{ j_1, \ldots, j_m \} , \; i_k = r_{j_k} , \; \beta_k = \frac{T}{N} \lambda_{j_k} \quad .$$

Then

$$s_{n,r_1,\ldots r_K}(q) = \binom{n}{i_1 \; . \; . \; . \; i_m} \left(1 - \sum_{k=1}^{m} \beta_k\right)^{n - \sum_{k=1}^{K} i_k} \prod_{k=1}^{m} \beta_k^{i_k} \quad (79)$$

The basic multinomial recurrence relation (see e.g. [13])

$$\binom{n+1}{i_1 \; . \; . \; . \; i_m} \quad (80)$$

$$= \frac{(n+1)!}{i_1!\cdots i_m!(n+1-\sum i_k)!} = \frac{n+1}{n+1-\sum i_k}\begin{pmatrix} & n & \\ i_1 & . \quad . \quad . & i_m \end{pmatrix}$$

$$= \begin{pmatrix} & n & \\ i_1 & . \quad . \quad . & i_m \end{pmatrix} + \sum_{k=1}^{m}\begin{pmatrix} & n & \\ i_1 & \dots \quad i_k - 1 \quad \dots & i_m \end{pmatrix}$$

shows that

$$\begin{pmatrix} & n+1 & \\ i_1 & . \quad . \quad . & i_m \end{pmatrix}\left(1-\sum_{k=1}^{m}\beta_k\right)^{n+1-\sum_{k=1}^{K} i_k}\prod_{k=1}^{m}\beta_k^{i_k}$$

$$= \left(1-\sum_{j=1}^{m}\beta_j\right)\begin{pmatrix} & n & \\ i_1 & . \quad . \quad . & i_m \end{pmatrix}\left(1-\sum_{k=1}^{m}\beta_k\right)^{n-\sum_{k=1}^{K} i_k}\prod_{k=1}^{m}\beta_k^{i_k}$$

$$+\sum_{j=1}^{m}\beta_j\begin{pmatrix} & n & \\ i_1 & \dots \quad i_j - 1 \quad \dots & i_m \end{pmatrix}\left(1-\sum_{k=1}^{m}\beta_k\right)^{n-\sum_{k=1}^{K} i_k+1}\beta_j^{i_k-1}\prod_{\substack{k=1\\k\neq j}}^{m}\beta_k^{i_k}$$

Thus

$$s_{n+1,r_1,\dots r_K}(q)$$
$$= s_{n,r_1,\dots,r_K}(q) - \frac{T}{N}\sum_{j\notin R_q}\lambda_j s_{n,r_1,\dots,r_K}(q) + \frac{T}{N}\sum_{j\notin R_q}\lambda_j s_{n,r_1,\dots,r_j-1,\dots,r_K}(q).$$

(81)

Summing over $A_{r,q}$ shows that the recursion equation is satisfied also for

$$s_{n+1,q_1,\dots q_K,r} = \sum_{\substack{r_1=\dots=r_K=0\\(r_1,\dots r_K)\in A_{r,q}}}^{n} s_{n+1,r_1,\dots r_K}(q)\pi_{n,q}$$

(82)

$$= \left(1-\frac{T}{N}\sum_{i\notin R_q}\lambda_i\right)s_{n,q_1,\dots q_K,r} + \frac{T}{N}\sum_{i\notin R_q}\lambda_i\sum_{\substack{r_1=\dots=r_K=0\\(r_1,\dots r_K)\in A_{r,q}}}^{n} s_{n,r_1,\dots r_i-1,\dots,r_K}(q)\pi_{n,q}$$

$$= \left(1-\frac{T}{N}\sum_{i\notin R_q}\lambda_i\right)s_{n,q_1,\dots q_K,r} + \frac{T}{N}\sum_{i\notin R_q}\lambda_i\sum_{\substack{r_1=\dots r_i'=\dots=r_K=0\\(r_1,\dots r_i',\dots r_K)\in A_{r,q-\omega_j}}}^{n} s_{n,r_1,\dots r_i',\dots,r_K}(q)\pi_{n,q}$$

$$= \left(1-\frac{T}{N}\sum_{i\notin R_q}\lambda_i\right)s_{n,q_1,\dots q_K,r} + \frac{T}{N}\sum_{i\notin R_q}\lambda_i s_{n,q_1,\dots q_K,r-\omega_i}$$

.

Let us sum over all values of $r$

$$\sum_{r=0}^{\infty} s_{n,q_1,\dots q_K,r} \tag{83}$$

$$= \sum_{r=0}^{\infty}\left(\sum_{\substack{r_1=\dots=r_K=0\\(r_1,\dots r_K)\in A_{r,q}}}^{n}\begin{pmatrix} & n & \\ r_1 & . \quad . \quad . & r_K \end{pmatrix}\left(1-\frac{T}{N}\sum_{i\notin R_q}\lambda_i\right)^{n-\sum_{i=1}^{K} r_i}\prod_{i\notin R_q}\left(\frac{T}{N}\lambda_i\right)^{r_i}\right)\pi_{n,q}$$

$$= \sum_{i_1=0}^{n}\cdots\sum_{i_m=0}^{n}\begin{pmatrix} & n & \\ i_1 & . \quad . \quad . & i_m \end{pmatrix}\left(1-\sum_{k=1}^{m}\beta_k\right)^{n-\sum_{k=1}^{K} i_k}\prod_{k=1}^{m}\beta_k^{i_k}\pi_{n,q}$$
$$= \pi_{n,q}.$$

Clearly $0 \le s_{n,q_1,\dots q_K,r} \le 1$. This completes the proof. □

*Proof of Theorem 3*: Write $\alpha_k = \frac{T}{N}\beta_k = \lambda_{j_k}$. The limit of the multinomial term

$$L = \lim_{N\to\infty}\begin{pmatrix} & n & \\ i_1 & . \quad . \quad . & i_m \end{pmatrix}\left(1-\sum_{k=1}^{m}\frac{T}{N}\alpha_k\right)^{xN-\sum_{k=1}^{K} i_k}\prod_{k=1}^{m}\left(\frac{T}{N}\alpha_k\right)^{i_k}$$

is calculated by Stirling's formula using Poisson's limit procedure:

$$L = \prod_{i\notin R_q}\frac{(xT\lambda_i)^{r_i}}{r_i!}e^{-xT\sum_{i\notin R_q}\lambda_i}. \tag{84}$$

Thus

$$s_{t,q_1,\dots q_K,r} = \lim_{N\to\infty} s_{xN,q_1,\dots,q_K,r} = \sum_{\substack{r_1=\dots=r_K=0\\(r_1,\dots r_K)\in A_{r,q}}}^{n}\prod_{i\notin R_q}\frac{(xT\lambda_i)^{r_i}}{r_i!}e^{-xT\sum_{i\notin R_q}\lambda_i}\pi_q$$

as claimed. □

*Proof of Lemma 2:* Directly summing from (31) is difficult, but as (31) is the limit of (30), it can be calculated as $O(\frac{1}{N})$ limit from (30). To get $O(\frac{1}{N})$ limit, we only can take the terms $r_i = 1$, $r_j = 0$ for $i \neq j$, and $r_1 = \dots = r_K = 0$ since there is the term $\left(\frac{T}{N}\right)^{r_i}$ in (30). The term $r_1 = \dots = r_K = 0$ disappears as $r = 0$ and for the other terms $r = \omega_i$. Consequently, from (30)

$$\sum_{r=0}^{\infty} rs_{n,q_1,\dots,q_K,r} = \sum_{i \notin R_q} \omega_i \lambda_i n \left(1 - \frac{T}{N} \sum_{j \notin R_q} \lambda_j \right)^{n-1} \frac{T}{N} \lambda_i \pi_q + O(N^{-2})$$

$$= \sum_{i \notin R_q} \omega_i \lambda_i n \frac{T}{N} \lambda_i \pi_q + O(N^{-2}) \cdot \qquad (85)$$

Thus

$$\sum_{r=0}^{\infty} rs_{t,q_1,\dots,q_K,r} = \lim_{N \to \infty} \sum_{r=0}^{\infty} rs_{xN,q_1,\dots,q_K,r}$$

$$= \lim_{N \to \infty} \sum_{i \notin R_q} \omega_i \lambda_i x T \pi_q = t \sum_{i \notin R_q} \omega_i \lambda_i \pi_q \cdot \square \qquad (86)$$

## ACKNOWLEDGMENT

Author is deeply grateful to Dr. Markku Liinaharja from Helsinki University of Technology (HUT) for reading the proofs & Dr. Jouni Karvo, from HUT, for suggesting improvements to the structure and readability of the paper, and to the anonymous reviewers for their helpful comments. Any possibly remaining errors are solely on author's responsibility,